\documentclass[review]{elsarticle}

\usepackage{lineno,hyperref}
\modulolinenumbers[5]

\usepackage{bigstrut}
\usepackage{amssymb,amsmath}
\usepackage{amsfonts}
\usepackage{multirow}
\usepackage{tabularx}
\usepackage{latexsym}
\usepackage{algorithm}
\usepackage{amsmath}
\usepackage{amssymb}
\usepackage{amsbsy}
\usepackage{float}
\usepackage{amsthm}
\usepackage{graphicx}
\usepackage{color}
\usepackage{float}
\usepackage{graphicx}
\usepackage{geometry}
\usepackage{algpseudocode}
\usepackage{listings}
\usepackage{xcolor}
\usepackage{algorithm}
\usepackage{algpseudocode}
\usepackage{tikz}
\usetikzlibrary{intersections}
\usepackage{tikz}
\usetikzlibrary{fadings}

\begin{document}

	\begin{frontmatter}
		
		\title{A Deterministic Model for the Transshipment Problem of a Fast Fashion Retailer under Capacity Constraints}
		
		
		\author[mymainaddress]{Siamak Naderi\corref{mycorrespondingauthor}}
		\cortext[mycorrespondingauthor]{Corresponding author}
		\ead{siamak@sabanciuniv.edu}
		\author[mymainaddress]{Kemal Kilic}
		\ead{kkilic@sabanciuniv.edu}
		\author[mysecondaryaddress]{Abdullah Dasci}
		\ead{dasci@sabanciuniv.edu}

		\address[mymainaddress]{Faculty of Engineering and Natural Sciences, Sabanci University, Orta Mahalle, Tuzla, Istanbul, Turkey}
		\address[mysecondaryaddress]{School of Management, Sabanci University, Orta Mahalle, Tuzla, Istanbul, Turkey }
		
		\begin{abstract}
			
			In this paper we present a novel transshipment problem for a large apparel retailer that operates an extensive retail network. Our problem is inspired by the logistics operations of a very large fast fashion retailer in Turkey, LC Waikiki, with over 450 retail branches and thousands of products. The purpose of transshipments is to rebalance stocks across the retail network to better match supply with demand. We formulate this problem as a large mixed integer linear program and develop a Lagrangian relaxation with a primal-dual approach to find upper bounds and a simulated annealing based metaheuristic to find promising solutions, both of which have proven to be quite effective. While our metaheuristic does not always produce better solutions than a commercial optimizer, it has consistently produced solutions with optimality gaps lower than 7\% while the commercial optimizer may produce very poor solutions with optimality gaps as high as almost 300\%. We have also conducted a set of numerical experiments to uncover implications of various operational practices of LC Waikiki on its system's performance and important managerial insights.   
			
		\end{abstract}
		
		
	\end{frontmatter}

	\section{Introduction}
	
	Due to its impact on revenues, costs, and more importantly, on service levels, logistics management has become increasingly critical in the apparel industry (Kiesmuller and Minner 2009\nocite{kiesmuller2009inventory}). As consumers demand greater product variety and higher levels of responsiveness at lower prices, effective management of logistics activities arises as a key competitive advantage for the retailers in this industry. The main challenges faced by these retailers are short selling seasons and unpredictable demands. Since forecasts are mostly inaccurate, firms usually have either excess inventories that are sold at markdown prices or stock-outs that lead to lost sales. The problem is exacerbated with short selling seasons which prevent firms to replenish their stocks. Therefore, an effective logistics strategy is key to avoid both of these undesirable outcomes.\\
	
	Logistics decisions of apparel retailers include initial ordering before the season begins, allocation to the branches at the beginning of the season, and eventually phasing-out of the products at the end of the selling season. Increasingly, however, retailers are also practicing what is called ``transshipment" or ``transfer" policies, which involve the reallocation of products among retail branches in mid-season (Li et al. 2013). These policies help retailers to reduce stock-outs as well as excess inventories. This is the issue that is addressed in this paper. \\
	
	The problem that we consider here is inspired by the logistics operations at the largest apparel retailer in Turkey, LC Waikiki, which has positioned itself as a ``fast fashion" retailer. The term fast fashion used to refer to inexpensive designs that appeared on catwalks and were quickly moved to store shelves. Fast fashion items based on the most recent trends have shaped mass-merchandized clothing collections. Therefore, mass-merchendize retailers compete to introduce latest fashion trends in their collections. Although the term was first used in the US in the 1980s, the expression did not receive worldwide adoption until popularized by the Spanish-based apparel giant Zara. The crucial issue in fast fashion is providing inexpensive collections that also respond to fast changing consumer tastes and trends. Therefore, the entire fast fashion supply chain must be sufficiently agile to operate with products for which life cycles are measured not in months but rather in weeks.\\
	
	On the plus side, the speed at which fast fashion moves tends to help retailers avoid markdowns. Typically, these retailers do not place very large orders months before the actual selling season, but rather work with smaller initial orders and renew collections more frequently. On the negative side, however, the fast-paced environment calls for higher turnover and more frequent introduction of new designs, a setting that necessitates shorter design and production lead times. As a result, companies need to rely on more expensive local sources and accommodate large design teams. This fast-paced environment also creates new logistics challenges for retailers: When will these products be replaced? Should they be completely removed from the stores or kept at display at select stores? What will happen to the leftover items; reintroduced elsewhere, sold at discount, or simply written-off? Fast fashion companies need to deal with these issues much more frequently than traditional retailers. \\
	
	Facing such challenges, leading fast fashion companies such as Zara and another Spanish company, Mango, the Japanese World Co., and Swedish H\&M have built supply chains aiming at quickly responding to consumers' changing demands while decreasing the excess inventories at branches and hence, lowering costs (Caro and Gallien 2007)\nocite{caro2007dynamic}. For instance, Zara developed a decision support system featuring demand updating and a dynamic optimization module for initial shipment decisions to avoid stock-outs as well as excess inventories (Gallien et al. 2015\nocite{gallien2015initial}). In addition to correct initial shipment decisions, the transfer or transshipment decisions among retail locations are also instrumental to reduce stock-outs and excess inventories.\\
	
	The main benefit of transfer actions is better matching inventory and demand at different locations. It uses up-to-date sales information and inventory data and redistribute available inventory among the retail locations. Due to socio-economic and geographical differences among retailer locations, it is possible that a product sells very well in some stores while less so in others. Transfer actions can be adopted as a tool to increase the inventory levels at receiver stores while providing extra shelf space at sender stores. This action can be adopted by bypassing the central depot to facilitate the quick movement of merchandise. As a result, the revenues are increased while costs are reduced as compared to a system where no transshipment is utilized (Tagaras 1989). There are a number works in the literature that describe how retailers take advantage of transfers to improve their performances. For example, Archibald et al. (2009)\nocite{archibald2009index} and Archibald et al. (2010)\nocite{archibald2010use} address transshipment issues at a tire retailer that has a network of 50 stores. In another work, Hu and Yu (2014) \nocite{hu2014optimization} present a proactive transshipment problem for a famous fashion brand in China that has network for 43 retailers in Shanghai. The problem that we introduce here is motivated by the largest apparel retailer in Turkey. \\
	
	In the next section, we provide a detailed background of our problem that includes the transfer practices at the company that motivated this work and a detailed description of the problem setting. In Section \ref{sec_LR} we give a brief literature review. Section \ref{PD} presents the progressive development of the mathematical model. Section \ref{SA} presents our solution methods that include a Lagrangian relaxation based upper bounding method and simulated annealing based metaheuristic to find good feasible solutions. Section \ref{CR} reports on our numerical experiments followed by a few concluding remarks in Section \ref{Concl}. \\
	
	\section{Background}\label{BG}
	
\begin{figure}[!b].
		\hspace*{-3em}
		\begin{minipage}{0.5\linewidth}
			\includegraphics[width = 1.05\textwidth]{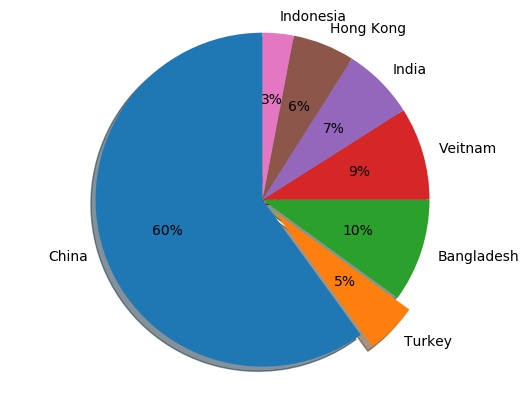}
		\end{minipage}
		\hfill
		\hspace*{-3em}
		\begin{minipage}{0.505\linewidth}
			\includegraphics[width = 1.05\textwidth]{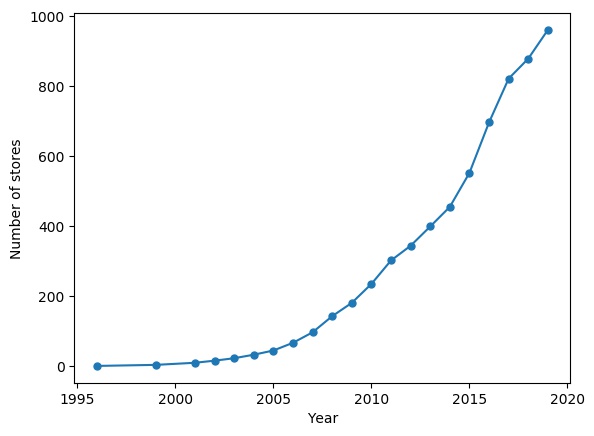}
		\end{minipage}
		\caption{Left: Global export market share, Right: Growth in the total number of stores}
		\label{growth_store}
	\end{figure}

	Textile is one of the key sectors in the Turkish economy in terms of GDP, domestic employment, and exports. Textile accounts for 10\% of the Turkish GDP and 20\% of employment in the manufacturing sector\footnote{blog.tcp.gov.tr}. In 2016 Turkey exported around 15 Billion USD, mainly to the European Union countries and was ranked as the 6th biggest textile exporting country (see Figure \ref{growth_store})\footnote{www.wikipedia.com}. LC Waikiki, which has provided motivation to this work, is the largest textile retailer in Turkey with significant international presence.  \\

	LC Waikiki was founded in 1988 in France by a French designer and his friend. The LC Waikiki brand name is created by adding the word Waikiki, a famous beach in Hawaii, to LC, the abbreviation of the French word ``Les Copains" meaning ``friends". TEMA, a Turkey based group which was then a major supplier of the company, bought the LC Waikiki brand in 1997 undertook a major restructuring that included focusing on domestic market. In the same year, the group entered the Turkish fashion retail market with 21 stores. In 2009, it opened its first international store in Romania since the TEMA group had purchased the brand. Over the years, the group has followed  an aggressive expansion strategy both domestically and internationally. Today, LC Waikiki has more than 370 stores in 34 countries in Asia, Africa, and Europe, in addition to over 480 stores in Turkey. In 2011, LC Waikiki became the leader of the ``Ready-to-Wear" market in Turkey and remains as the largest apparel retailer in terms of sales as well as the number of stores.  Figure \ref{growth_store} depicts LC Waikiki's phenomenal growth in terms of the total number of stores over the years. \\
	
	LC Waikiki has a highly centralized order planning and logistics system in which all initial orders and subsequent distribution decisions are made by the headquarters. New merchandise is received at a single central warehouse located in Istanbul, which then distributes essentially the entire amount to the retail branches (there are varying practices for international stores which are excluded from the consideration in this study). The retail practice at LC Waikiki can be considered as fast fashion in that it aims to keep items in stores only for about six to eight weeks. During this period, if the sales realize below expectations, they may reduce prices or if the sales display disparities across the stores, they may utilize transfers among the stores. Finally, at the end of their shelf-life, products are returned to the central warehouse and later sent to outlet stores (about 40 of the 480 stores are designated as outlet stores) or  simply given away to charities. Stock-outs and excess inventories are critical issues at LC Waikiki as in any fast fashion company due to forecast errors. Since LC Waikiki initially distributes all of the items to stores, transfer remains essentially as the only tool to deal with these issues by rebalancing inventories across the retail network. It is these transfer decisions that is the subject of this paper.\\
	
	Currently, a group at the headquarters manages transfer decisions. This group utilizes a mathematical model accompanied with some pre- and post-processing activities. However, we cannot disclose the precise nature of the model and the activities due to proprietary nature of these information. After transfer solutions are obtained, orders are automatically generated and transmitted to the stores. Store employees collect the products that have been chosen for transfer from the shelves and move them to a storage room. In the storage room, products are put in the boxes, each destined to a specific store without any re-assortments. Since the storage room capacities are limited, stores cannot to transfer more than what they can hold at their storage room. Once boxing is finished, the logistics company picks up the boxes and delivers them to their destinations. The boxes are ideally delivered before the weekend so that the transferred items can be put on shelves for the weekend sales.\\
	
	Although our work is motivated by LC Waikiki's logistics operations, we believe many of the features of our model would resonate with issues fast fashion retailers need to consider. In our model, we maximize a measure of the total profit which is the total revenue less the total logistics cost that includes transportation, handling, and inventory holding costs. We also include a number of operational constraints that represent the real practice of the company. For example, we consider a centrally managed system where stores have no control over the decisions, i.e., they may not refuse the transfer decisions. This is valid particularly for firms that own their stores and manage them centrally. We also restrict the total number of items and the total number of stores to which each store can make shipments. Both of these constraints are justified by the limited number of employees in the stores and sizes of the storage rooms. Furthermore, in our model once a product is decided to be transferred from one store to the other, the entire stock (all the available sizes) is sent to the same store. LC Waikiki justifies this practice by the simplicity of the picking operations, which otherwise would be too labor intensive. Here, we will also investigate the effects of these restrictions on system performance. \\
	
	There are a number of issues relevant to the fashion logistics decisions that we leave out of the scope of this work: i) Initial allocation decisions, ii) Uncertainty in demand, and iii) Dynamic nature of the decision making process. At LC Waikiki the initial shipment decisions are made after a pilot sales experiment in which they obtain sales information from about 30 stores. They then make initial allocations in which they essentially distribute the entire stock to the stores. Certainly, the option of transfers might impact the initial allocation. However, we believe that the impact is small due to two aspects in this logistics system. First, the company has a policy to allocate almost all of the available inventory to the stores keeping none at the central depot. Therefore, the firm cannot use central depot for reallocation of products. Second, since the company has flat transportation cost rate independent of origin-destination pair, regional risk pooling effect becomes irrelevant. Therefore, the impact of subsequent transfer practice on initial allocation decisions has lessened. Demand uncertainty is always a concern particularly in the fashion industry, and in fact, it is the demand uncertainty that makes the transfer problem relevant. However, at the time a product is considered for transfer, there is demand information for at least for a couple of weekends. Therefore, the company is able to make much more accurate demand forecasts after this initial sales information, as compared to the time the initial allocation decisions are made. Finally, the transfer problem ideally should consider the fact that transfer decisions are made every week and hence, there are subsequent recourse opportunities. However, considering a multi-stage decision environment under demand uncertainty is simply beyond analysis for the sizes that we envision, particularly with complicating operational constraints. Instead, we envision a setting where the firm makes demand forecast until a product is planned to stay on shelves and the transfer problem is solved on a rolling-horizon basis. This setting, we believe, is a reasonable compromise given the other complexities of the system.  Similarly, initial replenishment decisions are also important and they would be impacted by subsequent transfer options. However, considering transshipment and replenishment decisions jointly would also be extremely difficult considering the scale of our problem and particular operational constraints. As we will see, in the literature too, there are very few papers that consider these decisions jointly (see Paterson et al. 2011).  \\

	We have also assumed that each product's shelf life is known. This assumption is justified by the practice of the company where they keep merchandize for about six to eight weeks. Decision to continue displaying products on the shelves or removing them involves a number of other factors to consider. It requires information on new product designs as well as space considerations at the stores for different merchandize groups. These issues are also rather involved and therefore, kept out of the current study, but certainly worthwhile to consider in the future.

	\section{Related literature}\label{sec_LR}
	
	There is a vast literature as far back as the 1950's on lateral transshipment or, as we call here, transfer issues. Although both terms commonly describe the decisions considered here and we use them interchangeably, the term transshipment has a wider meaning and usage. Time and again, various studies have shown that transfer option between retailers improves supply chain performance in terms of costs, revenues, and service levels. For example, Tagaras (1989\nocite{tagaras1989effects}) shows that utilizing transfer in a system with two retail locations leads to significant cost reductions. Although transfers considerably increase transportation cost, systems with these options are superior to systems without them (Banerjee et al. 2003). \nocite{banerjee2003simulation} Furthermore, transfers enhance customer service levels without the burden of carrying extra safety stock at retail locations (Burton and Banerjee 2005)\nocite{burton2005cost}. \\
	
	There are essentially two types of transfers: emergency or reactive transfers and preventive or proactive transfers, which are differentiated mainly with respect to  their timing (Lee et al. 2007\nocite{lee2007effective}, Paterson et al. 2011\nocite{paterson2011inventory}, Seidscher and Minner 2013\nocite{seidscher2013semi}, and Ahmadi et al. 2016\nocite{ahmadi2016bi}). Reactive transfer refers to responding to realized stock-outs at a retail location by using available inventory at another location whereas proactive transfer refers to redistribution of inventories among locations before the actual demand is realized. The literature can be classified primarily along this dimension, although there are also works that consider them jointly. \\
	
	Perhaps the earliest work that considers reactive transfers is by Krishnan and Rao (1965) \nocite{krishnan1965inventory} who study a centralized one-echelon inventory system with the objective of minimizing the total cost through transfers. One of the main motivations for reactive transfer models comes from spare parts distribution systems for repairable items, as exemplified by one of the more notable earlier works by Lee (1987) \nocite{lee1987multi} who studies a single-echelon model, which is then extended by Axs{\"a}ter (1990\nocite{axsater1990modelling}) to a two-echelon system. More recent works on spare parts systems can be attributed to van Wijk et al. (2019\nocite{van2019optimal}) who consider a two-location system with lateral transshipment as well as an outside emergency option and Boucherie et al. (2018\nocite{boucherie2018two}) who consider a complex two-echelon inventory system with multiple local warehouses.\\
	
	Models with reactive transfer policies have also been studied for non-repairable items. A notable contribution is due to  Robinson (1990), \nocite{robinson1990optimal} who provides structural results for a two-retailer system and develops a heuristic for the initial ordering decisions considering the subsequent transshipments. Herer et al. (2006\nocite{herer2006multilocation}) extend this work by considering more general cost structures and {\"O}zdemir et al. (2013) \nocite{ozdemir2013multi} extend it considering capacity constraints on the transportation network. More recently, transshipment policies in systems with perishable items have also attracted research (see for example, Nakandala et al. 2017 \nocite{nakandala2017lateral} and Dehghani and Abbasi 2018\nocite{dehghani2018age} for such recent works).\\
	
	Proactive transfer is based on the concept of inventory rebalancing and is mostly utilized in periodic review inventory control framework. Allen (1958\nocite{allen1958redistribution}) provides perhaps the earliest model that considers proactive transfers in a single-period setting, which is then generalized by Das (1975) \nocite{das1975supply} who also considers the initial replenishment decision. There are also models that include the timing of the transshipment as decision in a dynamic setting (Agrawal et al. 2004\nocite{agrawal2004dynamic} and Tiacci and Saetta 2011\nocite{tiacci2011heuristic}) as well as in a static setting (Kiesmuller and Minner 2009\nocite{kiesmuller2009inventory}).\\
	
	Although the type of transfers may be dictated by operational conditions of the setting, proactive transfer policies are found to be superior to purely reactive policies both in terms of costs and stock-out levels (see for example, Banerjee \nocite{banerjee2003simulation} et al. 2003 and Burton and Banerjee 2005\nocite{burton2005cost}). In some settings, however, companies may also have opportunities to implement these policies jointly (see for example, Lee et al. 2007\nocite{lee2007effective} for such a model). Finally, although all of the works mentioned above and majority of research on the transshipment issues, assume that the systems are centrally operated decentralized systems where retailers might refuse transshipment requests have also attracted research recently (see for example, \c{C}\"{o}mez et al. 2012 \nocite{ccomez2012season} and Li et al. 2013\nocite{li2013coordinating}).\\
	
	As we have noted earlier, the literature on transshipment issues is vast with considerable growth in the last two decades. Since reviewing this voluminous literature is not possible here, we have only offered a very selective review. Aside from the types of the transfer (i.e., reactive vs. proactive), the literature on transshipment is also divided along two other important dimensions: Whether the models consider only transshipment decisions or jointly with replenishment decisions and whether they consider multiple locations or just two locations. We have classified aforementioned works and few others with respect to these characteristics as shown in Table \ref{literature}. We choose to put some classic and some more recent ones, but it is still  far from portraying a complete picture. We refer the reader to a somewhat older, but an excellent review by Paterson et al. (2011) \nocite{paterson2011inventory} who also provide a more thorough classification and a comprehensive review up to its publication date. \nocite{karmarkar1977one} \nocite{diks1996controlling} \nocite{herer1999lateral} \nocite{tagaras2002effectiveness} \nocite{nonaas2007optimal} \nocite{bendoly2004integrated}\\
	
	\begin{table}[!b]
		\caption{Main characteristics of reviewed transfer-related papers}
		\setlength{\tabcolsep}{.05em}
		\scriptsize{
			\begin{tabular}{ l l|m{3cm}|m{3.5cm}|m{2.8cm}|m{4.1cm}| }
				\cline{3-6}
				& & \multicolumn{2}{c }{\textbf{Single Period}} & \multicolumn{2}{ |c| }{\textbf{Multiple Period}}\\
				\cline{3-6}
				& \textbf{Replenishment}&\multicolumn{1}{c|}{\textbf{2 Retailers}}  & \multicolumn{1}{c|}{\textbf{Multiple Retailers}} & \multicolumn{1}{c|}{\textbf{2 Retailers}} & \multicolumn{1}{c|}{\textbf{Multiple Retailers}} \\ \hline
				\multicolumn{1}{|c}{\multirow{2}{*}{\textbf{Proactive}}}& \multicolumn{1}{|c|}{\textbf{Yes}} &   Das 1975\newline Tagaras and Vlachos 2002  & Karmarkar and Patel 1977 \newline Hoadley and Heyman 1977\nocite{hoadley1977two} & Tiacci and Saetta 2011 \newline Abouee-Mehrizi et al. 2015\nocite{abouee2015optimal}  & Diks and Kok 1996 \newline  Ahmadi et al. 2016 \newline Feng et al. 2017\nocite{feng2017preventive} \\ \cline{2-6}
				\multicolumn{1}{|c}{} & \multicolumn{1}{|c|}{\textbf{No}} & Kiesmuller and Minner 2009 \newline Li et al. 2013  & Allen 1958 \newline Agrawal et al. 2004 &  Dan et al. 2016\nocite{dan2016ordering} &   Bertrand and Bookbinder 1998\nocite{bertrand1998stock} \newline Banerjee et al 2003 \newline Burton and Banerjee 2005 \newline Acimovic and Graves 2014\nocite{acimovic2014making}\newline Peres et al. 2017 \\ \cline{1-6}
				\multicolumn{1}{|c}{\multirow{2}{*}{\textbf{Reactive}}}& \multicolumn{1}{|c|}{\textbf{Yes}} & Herer and Rashit 1999 \newline Minner and Silver 2005\nocite{minner2005evaluation} \newline Liao et al. 2014\nocite{liao2014optimal}\newline Olsson 2015\nocite{olsson2015emergency} \newline Dehghani and Abbasi 2018 & Lee 1987 \newline Axs{\"a}ter 1990 \newline Herer et al. 2006 \newline Johansson and Olsson 2018\nocite{johansson2018age} \newline Boucherie et al. 2018 & Archibald et al. 1997\nocite{archibald1997optimal} \newline Herer and Tzur 2001\nocite{herer2001dynamic} \newline van Wijk et al. 2019\nocite{van2019optimal} &   Archibald et al. 2009 \newline van Wijk et al. 2012\nocite{van2012approximate} \newline {\"O}zdemir et al. 2013\\ \cline{2-6}
				\multicolumn{1}{|c}{} & \multicolumn{1}{|c|}{\textbf{No}} & Herer and Rashit 1995\nocite{herer1995lateral} \newline Shao et al. 2011\nocite{shao2011incentives} \newline Liao et al. 2014\nocite{liao2014optimal} & Non{\aa}s and J{\"o}rnsten 2007  \newline Hu and Yu 2014  \newline Patriarca et al. 2016\nocite{patriarca2016inventory} \newline Bhatnagar and Lin 2019\nocite{bhatnagar2019joint}&Tagaras 1989 \newline Comez et al. 2012 \newline Shao 2018\nocite{shao2018production}  & Robinson 1990 \newline Banerjee et al. 2003 \newline Burton and Banerjee 2005 \newline Dijkstra et al. 2017\nocite{dijkstra2017transshipments}\\ \cline{1-6}
			\end{tabular}}
			\label{literature}
		\end{table}

		Despite many simplification attempts, solving transfer problem to optimality remains a challenge. Even in the presence of many simplifying assumptions such as single product, single-period, limited number of retail locations, static transfer timing and so on, past works can only provide approximate solutions. The problem that we present here considers proactive transfers, but since it is motivated by the actual logistics operations at a large fashion retailer, it has many complexities that would be quite challenging to resolve under demand uncertainty or in a dynamic fashion. Therefore, we need to make restrictive assumptions along these dimensions. Certainly we are not alone in this respect; there are numerous other works that consider deterministic demand for transshipment models. Not surprisingly, these works also contain other complicating factors. For example, Herer and Tzur (2001\nocite{herer2001dynamic} and 2003\nocite{herer2003optimal}) in their multi-period model consider fixed ordering costs in transshipments; Lim et al. (2005\nocite{lim2005transshipment}) and Ma et al. (2011)\nocite{ma2011crossdocking} study transshipment decisions via cross-docking locations under time windows; Qi (2006\nocite{qi2006logistics}) considers transshipment and production scheduling decisions jointly; Lee (2015\nocite{lee2015modeling}) considers concave production and transportation costs; Coelho et al. \nocite{coelho2012inventory} (2012), Mirzapour Al-e-hashem and Rekik (2014)\nocite{mirzapour2014multi}, and Peres et al. (2017\nocite{peres2017optimization}) consider routing issues alongside transshipments; Rahmouni et al. (2015\nocite{rahmouni2015multi}) and Feng et al. (2017) develop EOQ-based delivery scheduling models with transshipment while considering multiple products and resource constraints. Our setting too has a few operational practices that force us to model a static and deterministic problem.

		\section{Problem description and model formulation}\label{PD}
		
		We consider a retail logistics system that consists of a number of retail stores, each carrying a set of products of different sizes (SKUs). The firm has the precise stock information; that is, how many of each SKU the stores have and the projected demands of each SKU at each location during the remainder of the sales period. The problem is how to reallocate (some of) the products to maximize a profit measure that is total revenue less transfer, handling, and inventory holding costs. \\
		
		The firm has a single price policy in that the same price is applied to a product at all locations, which is indeed the practice of many retail chains and particularly of LC Waikiki. Each product also has a fixed transfer cost regardless of the origin-destination pair. This assumption is also motivated by the practice at LC Waikiki which has outsourced the transportation operations to a logistics company. The transfers are made by standard sized boxes for which LC Waikiki pays a fixed amount regardless of its contents and the locations of the sender and receiver stores. Since the number of products that fit in a box depends on the volume of the product, transportation cost differs for each product, but not on origin-destination pairs. The transfer cost can be estimated by adding the handling cost to the transportation cost for each product. However, none of these assumptions are really essential either for modeling or for our solution method and they can easily be relaxed.\\
		
		We also assume that transfer time has no effect on the sales. The main purpose is to finalize the delivery of transfer items before the weekend where the bulk of the sales materialize. Hence, delivering a day earlier or later presumably does not make much difference, as long as the products arrive for the weekend. Furthermore, the geography of Turkey does not allow wide variations in transfer times, but we also recognize that considering transfer time effects would be valuable in some settings. Finally, we assume that there are no replenishment opportunities from the central warehouse at the time of the transfer decisions. Since the company has a policy to allocate the entire inventory of a product to the stores at the beginning rather than keeping some at the warehouse, this assumption is well justified. As another operational practice, they do not consider a second replenishment option. This is a common practice among the fast-fashion retailers whose business practices involve speedy turnover of designs as exemplified by Zara's practice (see for example, Gallien et al. 2015). \\
		
		In addition to these requirements, we assume a single-period setting and deterministic demand. At LC Waikiki, most of the sales occur at weekends and therefore, the inventory levels of each product are updated at the beginning of each week. Likewise, demand forecasts are also revised after observing weekend sales. As a result, LC Waikiki, solves the transshipment problem once a week which allows us to consider single-period assumption to decrease the complexity of the problem. Deterministic demand is assumed since the forecast from the company is fairly accurate. After two or three weekend sales, the company can have a fairly good idea about the demand in the rest of the products' shelf lives. It is stated that their forecast error is below 15\%. Many papers related to fast-fashion also state that forecast errors are considerably smaller towards the end of shelf lives of products (see for example, Caro and Gallien, 2010\nocite{caro2010inventory}) . Finally, as mentioned earlier, the company has a few operational practices that we include in our model: If a product is transferred from a store, its entire available inventory (all SKUs) is shipped to a single store. Also, there are limits on the total number of SKUs that can be transferred from a store and the number of different destinations to which a store can make transfers. All these assumptions could be relaxed or generalized, but we choose to stay with the company practices as much as possible.\\
		
		As we will see shortly, without the aforementioned operational constraints, the problem can simply be formulated as a profit-maximizing transportation problem, which can easily be solved as a linear program. We are also ensured integer solutions if the demand and inventory values are integers. When we add the restriction on the total number of SKUs that can be transferred from a store, the problem can still be solved as a linear program. When we further add the restriction on the number of stores that a store can ship to, however, we need to introduce binary variables to keep track of whether a shipment is made from one store to another. Finally, when we include the single-destination constraint (i.e., when a product is shipped from one store to another, all the SKUs of the product must be shipped) the problem becomes much more difficult because we now also need to define a much larger set of binary decision variables to keep track of shipments between stores.\\
		
		We now give the preliminary definitions, followed by the formulation of the model. We start with the base model without considering the operational requirements of the company and progressively extend the model by adding each of these constraints. We call two stores as ``connected" if at least one product is transferred from one store to the other.
		
		\subsection*{Sets and indices:}
		
		\begin{itemize}
			\item[] $i,j \in I:$ Set of stores,
			\item[] $p \in P:$ Set of products,
			\item[] $k \in K_{p}:$ Set of sizes for each product $p \in P$.
		\end{itemize}
		
		\subsection*{Parameters:}
		
		\begin{itemize}
			\item[] $s_{ipk}:$ Stock level of size $k$ of product $p$ at store $i$,
			\item[] $d_{ipk}:$ Demand of size $k$ of product $p$ at store $i$,
			\item[] $r_p$: Unit net revenue of product $p$,
			\item[] $c_p$: Unit transfer cost of product $p$,
			\item[] $h_{p}:$ Holding cost of product $p$.
		\end{itemize}
		
		\subsection*{Decision variables:}
		
		\begin{itemize}
			\item[] $x_{ijpk} :  \mbox{ Amount of size } k \mbox{ of product } p \mbox{ transfered from store  } i \mbox{ to store } j,$
			\item[] $z_{ipk} :  \mbox{ Sales of size } k \mbox{ of product } p \mbox{ at store } i,$
			\item[] $w_{ipk} :  \mbox{ Amount of size } k \mbox{ of product } p \mbox{ store } i \mbox{ has after the transfers}.$
		\end{itemize}

		\subsection*{Relaxed model:}
		\begin{subequations}
			\begin{eqnarray}
				\label{relaxed_objective}\max \Pi & = & \sum_{j \in I}\sum_{k \in k_{p}} \sum_{p \in P}r_{p} z_{jpk} - \sum_{i \in I} \sum_{\substack{j \in I \\ j\neq i}} \sum_{p \in P}\sum_{k \in K_p} c_p x_{ijpk} - \sum_{i \in I } \sum_{p \in P } \sum_{k \in K_p }h_{p}(w_{ipk}-z_{ipk}) \\
				\label{relaxed_recieved} \mbox{s.t.} & & w_{ipk} = \sum_{j \in I}x_{jipk}, \mbox{ for all } i \in I , p \in P \mbox{ and } k \in K_{p},\\
				\label{relaxed_sales_supply} & & z_{ipk} \leq w_{ipk} , \mbox{ for all } i \in I, p \in P \mbox{ and } k \in K_p, \\
				\label{relaxed_sales_demand} & & z_{ipk} \leq d_{ipk}, \mbox{ for all } i \in I, p \in P \mbox{ and } k \in K_p, \\
				\label{relaxed_shipment_number2} & &  x_{ijpk} \leq s_{ipk} , \mbox{ for all } i,j \in I, p \in P \mbox{ and } k \in K_p,  \\
				\label{relaxed_binary_x} & & x_{ijpk} \geq 0, \mbox{ for all } i,j \in I, p \in P\mbox{ and } k \in K_p  ,\\
				\label{relaxed_binary_z} & & z_{ipk} \geq 0, \mbox{ for all } i \in I, p \in P \mbox{ and } k \in K_p.
			\end{eqnarray}
		\end{subequations}
		
		The objective function (\ref{relaxed_objective}) maximizes the total profit where the first term represents the total revenue obtained from sales, the second term is the total transfer cost, and the last term is the total holding cost. Constraints (\ref{relaxed_recieved}) define the stock level of each SKU after the transfers are completed. Constraints (\ref{relaxed_sales_supply}) and (\ref{relaxed_sales_demand}) ensure that sales are less than or equal to demand or the available stock of SKUs after the transfers are made. Constraints (\ref{relaxed_shipment_number2}) guarantee that a store may not transfer more than its inventory. Constraints (\ref{relaxed_binary_x}) and (\ref{relaxed_binary_z}) define the decision variables.\\

		As mentioned earlier, above problem is simply a profit maximizing transportation problem and can be easily solved by commercial optimizers. Now we extend the problem (\ref{relaxed_objective})-(\ref{relaxed_binary_z}) by adding one of the capacity constraints:

		\begin{subequations}
			\begin{eqnarray}
				\label{relaxed_objective_ex1}\max &(\ref{relaxed_objective})& \\
				\label{relaxed_knapsack_ex1} \mbox{s.t.} & & \sum_{\substack{j \in I \\ j\neq i}}\sum_{p \in P}\sum_{k \in K_p} x_{ijpk} \leq  A_i, \mbox{ for all } i \in I, \\
				\label{all_const_ex1} & & (\ref{relaxed_recieved})-(\ref{relaxed_binary_z}).
			\end{eqnarray}
		\end{subequations}
		where Constraints (\ref{relaxed_knapsack_ex1}) ensure that a store does not transfer more SKUs than it is allowed. This constraint does not pose a challenge as the problem is still a linear program. \\
		
		Next, we add the second capacity constraint to the current model. To do so, however, we need to introduce a binary decision variable $y_{ij}$ that represents if stores $i$ and $j$ are connected. The extended model is formulated as follows:
		
		\begin{subequations}
			\begin{eqnarray}
				\label{relaxed_objective_ex2}\max &(\ref{relaxed_objective})& \\
				\label{relaxed_knapsack_ex2} \mbox{s.t.} & & \sum_{\substack{j \in I \\ j\neq i}}\sum_{p \in P}\sum_{k \in K_p} x_{ijpk} \leq  A_i, \mbox{ for all } i \in I, \\
				\label{relaxed_shipment_number2_ex2} & & \sum_{\substack{j \in I \\ j\neq i}} y_{ij} \leq B_i , \mbox{ for all } i \in I, \\
				\label{basic_shipment_number1_ex2} & & x_{ijpk} \leq s_{ipk} y_{ij} , \mbox{ for all } i,j \in I, j\neq i,p \in P \mbox{ and } k \in K_p  , \\
				\label{all_const_ex2} & & (\ref{relaxed_recieved})-(\ref{relaxed_sales_demand}),\; (\ref{relaxed_binary_x}), \mbox{ and } (\ref{relaxed_binary_z}).
			\end{eqnarray}
		\end{subequations}
		where Constraints (\ref{relaxed_shipment_number2_ex2}) do not allow a particular store to transfer to more than a given number of stores. Constraints (\ref{basic_shipment_number1_ex2}) allow transfer between two stores only if they are connected; these constraints essentially replace Constraints (\ref{relaxed_shipment_number2}).\\
		
		Finally, single-destination constraint is added to the model. This constraint requires a change in one decision variable set that represents the SKU flow. We now define a binary decision variable $x_{ijp}$ that represents if product $p$ is transfered from store $i$ to store $j$, or not. Then, $x_{ijpk} = s_{ipk} x_{ijp}$, which allows us to drop the original flow variables from the formulation. The final model is given below. 
		
		\subsection*{The final model:}
		\begin{subequations}
			\small
			\begin{eqnarray}
				\label{basic_objective_nl}\max \Pi & = & \sum_{i \in I}\sum_{k \in k_{p}} \sum_{p \in P}r_{p} z_{ipk} - \sum_{i \in I} \sum_{\substack{j \in I\\ j\neq i}} \sum_{p \in P} \sum_{k \in K_p} c_p s_{ipk} x_{ijp} - \sum_{i \in I } \sum_{p \in P } \sum_{k \in K_p }h_{p}(w_{ipk}-z_{ipk}) \\
				\label{basic_sales_supply} \mbox{s.t.} & & z_{ipk} \leq \sum_{j \in I} s_{jpk} x_{jip}, \mbox{ for all } i \in I, p \in P \mbox{ and } k \in K_p, \\
				\label{basic_sales_demand} & & z_{ipk} \leq d_{ipk}, \mbox{ for all } i \in I, p \in P \mbox{ and } k \in K_p, \\
				\label{PM_pw_sales}  & & w_{ipk} = \sum_{j \in I}x_{jip}s_{jpk}, \mbox{ for all } i \in I , p \in P \mbox{ and } k \in K_{p},\\
				\label{basic_assignment1}  & & \sum_{j \in J} x_{ijp} = 1, \mbox{ for all } i \in I, p \in P, \\
				\label{basic_knapsack} & & \sum_{\substack{j \in I \\ j \neq i}} \sum_{p \in P} \sum_{k \in K_p} s_{ipk} x_{ijp} \leq A_i, \mbox{ for all } i \in I, \\
				\label{basic_shipment_number2} & & \sum_{\substack{j \in I \\ j \neq i}} y_{ij} \leq B_i , \mbox{ for all } i \in I, \\
				\label{basic_shipment_number1} & & x_{ijp} \leq y_{ij} , \mbox{ for all } i,j \in I, j\neq i \mbox{ and }  p \in P, \\
				\label{basic_binary_x} & & x_{ijp} \in \{0,1\}, \mbox{ for all } i,j \in I \mbox{ and } p \in P,\\
				\label{basic_binary_y} & & y_{ij} \in \{0,1\}, \mbox{ for all } i,j \in I,\\
				\label{basic_nonzero_z} & & z_{ipk} \geq 0, \mbox{ for all } i \in I, p \in P \mbox{ and } k \in K_p.
			\end{eqnarray}
		\end{subequations}
		
		 where Constraints (\ref{basic_assignment1}) ensure that if a product is transferred from a store, its entire inventory is moved to exactly one store. As a result, assignment to multiple stores is not allowed and similarly, a store may not also keep a portion of the inventory. As we have elaborated before, this ``single-destination'' practice is rather peculiar, but nonetheless it is the case at LC Waikiki.  The company justify this practice on the grounds that without this they would have to devote too much of their sales personnels' times for collection, which they are not willing to do. Clearly, this assumption may have a substantial impact on the profit, as it may severely restrict options to better match demand with the supply. Indeed, in our numerical experiments we try to give a sense of the implications of this assumption. As we have also noted, this assumption also complicates the problem substantially, without which the problem can be solved much more effectively.  \\
		
		Before we move to the analysis of the problem, we like to point out that the final model is indeed quite difficult. The following proposition shows that the problem is NP-hard.\\
		
		\newtheorem{prop}{Proposition}
		\begin{prop}
			Problem (\ref{basic_objective_nl})-(\ref{basic_nonzero_z}) is NP-hard.
		\end{prop}
		
		\emph{Proof:} We will prove the proposition by reduction. Assume that there is only one product ($P=\{1\}$),  no holding cost, ($h=0$) and the product has only one size ($K_{1}=\{1\})$. Furthermore, assume that the unit net revenue of the product is zero, ($r=0$), and there is no limitation on the number of stores to which each store can be connected (unlimited $B_{i}$). Since $r=0$, Constraints (\ref{basic_sales_supply}) and (\ref{basic_sales_demand}) become redundant. Moreover, if $B_{i}$ is unlimited, Constraints (\ref{basic_shipment_number2}) become redundant. Consequently, since any $y_{ij}$ can be one, Constraints (\ref{basic_shipment_number1}) are also redundant. Now the problem reduces to:
		\begin{proof}
		\begin{subequations}
			\begin{eqnarray}
				\label{gen_assg_obj} \min \Phi & = & \sum_{\substack{i,j \in I\\j\neq i}} c x_{ij}	\\
				\label{gen_assg_cons1} \mbox{s.t.} & & \sum_{j \in I} x_{ij}=1, \mbox{ for all } i \in I, \\
				\label{gen_assg_cons2} & & \sum_{\substack{j \in I\\j\neq i}}s_{i}x_{ij} \leq A_{i}, \mbox{ for all } i \in I ,\\
				\label{gen_assg_cons3} & & x_{ij} \in \{0,1\}, \mbox{ for all } i,j \in I.
			\end{eqnarray}
		\end{subequations}
		Problem (\ref{gen_assg_obj})-(\ref{gen_assg_cons3}) is the well-known generalized assignment problem which belongs to class of NP-hard problems (Savelsberg 1997\nocite{savelsbergh1997branch}).
		\end{proof}
		As the proposition shows, our problem (\ref{basic_objective_nl})-(\ref{basic_nonzero_z}) is a very difficult mixed integer linear problem. As we will present later, our experiments with a commercial optimizer demonstrated that this problem could not be solved effectively. At LC Waikiki, the number of products that are considered for transfer is about 2,000,  on average. On the other hand, the number of stores is approximately 450 nationwide. Therefore, the proposed mixed integer linear program can be huge and a heuristic approach seems to be a reasonable way to proceed. The next section describes such a heuristic method.\\

		\section{Solution approach}\label{SA}

		We have developed a Lagrangian Relaxation (LR) based approach to obtain good upper bounds in reasonable time. LR has shown exceptional success in solving large scale combinatorial optimization problems (Fisher 1981\nocite{fisher1981lagrangian}). LR is also used in the context of transshipment and it is shown that it can provide acceptable bounds to the optimal solution (Wong et al. 2005\nocite{wong2005simple} and Wong et al. 2006\nocite{wong2006multi}). A solution of the Lagrangian dual provides an upper bound on the optimal solution of the problem (\ref{basic_objective_nl})-(\ref{basic_nonzero_z}). To obtain a lower bound (i.e., a feasible solution), we have developed a two-stage heuristic that consists of a construction heuristic and simulated annealing based metaheuristic to improve the solution. Different heuristic and metaheuristic methods are applied to transshipment problems. For example, Patriarca et al. (2016) and Peres et al. (2017) develop metaheuristics to solve transshipment in inventory-routing problems. The latter applied a variable neighborhood search based algorithm, while the former developed a genetic algorithm. Moreover, local search based methods are utilized in transshipment problems. For instance, Wong et al. (2005) and Wong et al. (2006) developed a simulated annealing based metaheuristic to find promising feasible solutions. Therefore, we have also opted for such metaheuristic. In the rest of this section, we describe these methods in detail.\\

		\subsection{Obtaining upper bounds}
		
		Note that in the formulation, Constraints (\ref{basic_sales_demand}), (\ref{basic_shipment_number2}), and (\ref{basic_knapsack}) are similar to knapsack constraints and Constraints (\ref{basic_assignment1}) are basic assignment constraints, all of which are well-known in the literature. On the other hand, Constraints (\ref{basic_sales_supply}) and (\ref{basic_shipment_number1}) complicate the problem because they connect ``$z$" variables to ``$x$" variables and ``$x$" variables to ``$y$" variables, respectively. Thus, problem (\ref{basic_objective_nl})-(\ref{basic_nonzero_z}) can be decomposed in well-known problems by relaxing these complicating constraints. \\
		
		Let $\overline{\boldsymbol{\alpha}} = \{\alpha_{ipk} \in \mathbb{R}^{+}:i\in I, p \in P, k \in K_{p} \}$ and $\overline{\boldsymbol{\beta}}=\{\beta_{ijp}\in \mathbb{R}^{+}: i,j \in I,p \in P \}$ represent vectors of Lagrangian multipliers associated with Constraints (\ref{basic_sales_supply}) and (\ref{basic_shipment_number1}), respectively. Then the relaxed problem can be written as \\
		
		\begin{subequations}
			\small{
				\begin{eqnarray}
					\label{basic_objective_LR2}\max \Pi_{LR}(\overline{\boldsymbol{\alpha}},\overline{\boldsymbol{\beta}}) & = &  \sum_{i \in I} \sum_{p \in P} \sum_{k \in K_p} r_p z_{ipk} - \sum_{i \in I} \sum_{\substack{j \in I\\ j\neq i}} \sum_{p \in P}\sum_{k \in K_p} c_p s_{ipk} x_{ijp} - \sum_{i \in I } \sum_{p \in P } \sum_{k \in K_p }h_{p}(w_{ipk}-z_{ipk}) \notag \\
					& & -\sum_{i \in I}\sum_{p \in P}\sum_{k \in K_{p}} \alpha_{ipk}(z_{ipk}-\sum_{j \in I}s_{jpk}x_{jip}) - \sum_{i \in I}\sum_{\substack{j \in I \\ j\neq I}}\sum_{p \in P}\beta_{ijp}(x_{ijp}-y_{ij}) \\
					\label{Constraints} \mbox{s.t.} & & (\ref{basic_sales_demand})-(\ref{basic_shipment_number2}) \mbox{ and } (\ref{basic_binary_x})-(\ref{basic_nonzero_z}) .
				\end{eqnarray}}
			\end{subequations}
			
			\indent This problem can be decomposed into three subproblems, which are given as follows:
			\begin{subequations}
				\begin{eqnarray}
					\hspace*{-1.7cm}  \label{basic_objective_LR2_z} \mbox{\emph{Subproblem 1}: } \max \Pi_{LR}^{z}(\overline{\boldsymbol{\alpha}}) & = &  \sum_{i \in I} \sum_{p \in P} \sum_{k \in K_p} z_{ipk}(r_{p}-\alpha_{ipk}+h_p) \\
					\label{basic_sales_demand_LR2_z} \mbox{s.t.} & & z_{ipk} \leq d_{ipk}, \mbox{ for all } i \in I, p \in P \mbox{ and } k \in K_p, \\
					\label{basic_nonzero_z_LR2_z} & & z_{ipk} \geq 0, \mbox{ for all } i \in I, p \in P \mbox{ and } k \in K_p.
				\end{eqnarray}
			\end{subequations}
			
			\begin{subequations}
				\begin{eqnarray}
					\label{basic_objective_LR2_x} \mbox{\emph{Subproblem 2}: } \max \Pi_{LR}^{x}(\overline{\boldsymbol{\alpha}},\overline{\boldsymbol{\beta}}) & = & \sum_{i \in I}\sum_{\substack{j \in I\\ j \neq i}} \sum_{p \in P}x_{ijp}(\sum_{k \in K_{p}}((-c_{p}+\alpha_{jpk}-h_p)s_{ipk})-\beta_{ijp}) \notag \\
					& & +\sum_{i \in I}\sum_{p \in P}\sum_{k \in K_{p}}(\alpha_{iok}-h_p)s_{iok}x_{iip} + \sum_{i \in I}\sum_{j \in I} \sum_{p \in P} h_{p}w_{ipk}\\
					\label{basic_assignment1_LR2_x} \mbox{s.t.} & & \sum_{j \in J} x_{ijp} = 1, \mbox{ for all } i \in I, p \in P, \\
					\label{basic_knapsack_LR2_x} & & \sum_{\substack{j \in I \\ j \neq i}} \sum_{p \in P} \sum_{k \in K_p} s_{ipk} x_{ijp} \leq A_i, \mbox{ for all } i \in I, \\
					\label{PM_pw_sales_LR2_x}  & & w_{ipk} = \sum_{j \in I}x_{jip}s_{jpk}, \mbox{ for all } i \in I , p \in P \mbox{ and } k \in K_{p},\\
					\label{basic_binary_x_LR2_x} & & x_{ijp} \in \{0,1\}, \mbox{ for all } i,j \in I \mbox{ and } p \in P.
				\end{eqnarray}
			\end{subequations}
			
			\begin{subequations}
				
				\begin{eqnarray}
					\hspace*{-4.4cm}	\label{basic_objective_LR2_y} \mbox{\emph{Subproblem 3}: } \max \Pi_{LR}^{y}(\overline{\boldsymbol{\beta}}) & = &  \sum_{i \in I}\sum_{\substack{j \in I\\j\neq I}}\sum_{p \in P}\beta_{ijp}y_{ij} \\
					\label{basic_shipment_number2_LR2_y}\mbox{s.t.} & & \sum_{\substack{j \in I\\ j \neq i}} y_{ij} \leq B_i , \mbox{ for all } i \in I, \\
					\label{basic_binary_y_LR2_y} & & y_{ij} \in \{0,1\} \mbox{ for all } i,j \in I.
				\end{eqnarray}
			\end{subequations}

			\indent Among these problems, Problem (\ref{basic_objective_LR2_z})-(\ref{basic_nonzero_z_LR2_z}) is solvable by inspection. Problem (\ref{basic_objective_LR2_y})-(\ref{basic_binary_y_LR2_y}) can be decomposed into knapsack problems for each store. Problem (\ref{basic_objective_LR2_x})-(\ref{basic_binary_x_LR2_x}) seems to be computationally the most challenging of the three since this problem is similar to the generalized assignment problem. However, it is also separable for each store, which allows us to efficiently solve it. The subproblems for each $i \in I$ can be written as \\
			\begin{subequations}
				\begin{eqnarray}
					\label{basic_objective_LR3_x}\max \Pi_{LR}^{x^{i}}(\overline{\boldsymbol{\alpha}},\overline{\boldsymbol{\beta}}) & = & \sum_{\substack{j \in I \\ j \neq i}} \sum_{p \in P}x_{ijp}((\sum_{k \in K_{p}}(-c_{p}+\alpha_{jpk}-h_p)s_{ipk})-\beta_{ijp}) \notag \\
					& &  +\sum_{p \in P}\sum_{k \in K_{p}}(\alpha_{ipk}-h_p)s_{ipk}x_{iip} + \sum_{j \in I} \sum_{p \in P} h_{p}w_{ipk}\\
					\label{basic_assignment1_LR3_x} \mbox{s.t.} & & \sum_{j \in J} x_{ijp} = 1, \mbox{ for all } p \in P, \\
					\label{basic_knapsack_LR3_x} & & \sum_{\substack{j \in I\\ j \neq i}} \sum_{p \in P} s_{ipk} x_{ijp} \leq A_i, \\
					\label{PM_pw_sales_LR3_x}  & & w_{ipk} = \sum_{j \in I}x_{jip}s_{jpk}, \mbox{ for all }  p \in P \mbox{ and } k \in K_{p},\\
					\label{basic_binary_x_LR3_x} & & x_{ijp} \in \{0,1\}, \mbox{ for all } j \in I \mbox{ and } p \in P.
				\end{eqnarray}
			\end{subequations}

			Suppose that Lagrangian multipliers $\alpha_{ipk}$ and $\beta_{ijp}$ are set to some values. Then, let us define $\widehat{\text{z}}=\{\widehat{z}_{ipk}: i \in I, p \in P, k\in K_{p} \}$, $\widehat{\text{x}}=\{\widehat{x}_{ijp}: i,j \in I, p \in P \}$, and $\widehat{\text{y}}=\{\widehat{y}_{ij}: i,j\in I \}$ as the corresponding optimal solutions to the subproblems (\ref{basic_objective_LR2_z})-(\ref{basic_nonzero_z_LR2_z}), (\ref{basic_objective_LR2_x})-(\ref{basic_binary_x_LR2_x}) and (\ref{basic_objective_LR2_y})-(\ref{basic_binary_y_LR2_y}), respectively. We can then improve the Lagrangian bounds for a given solution by revising these Lagrangian multipliers. We achieve this by solving the Lagrangian dual while retaining the primal solutions. Interested reader can refer to Litvinchev (2007\nocite{litvinchev2007refinement}) for a detailed account of this approach. The Lagrangian dual can be formulated as follows:  \\

			\begin{subequations}
				\begin{eqnarray}
					\label{multipliers_objective}
					\min_{\overline{\alpha},\overline{\beta}} \max \Delta(\widehat{\text{z}},\widehat{\text{x}},\widehat{\text{y}})  & = & \sum_{i \in I} \sum_{p \in P} \sum_{k \in K_p} \widehat{z}_{ipk}(r_{p}-\alpha_{ipk}+h_p)+ \sum_{i \in I}\sum_{\substack{j \in I\\j\neq I}}\sum_{p \in P}\beta_{ijp}\widehat{y}_{ij} +\sum_{i \in I}\sum_{p \in P}\sum_{k \in K_{p}}(\alpha_{ipk}-h_p)s_{ipk}\widehat{x}_{iip} \notag \\
					& &   +\sum_{i \in I}\sum_{\substack{j \in I\\ j \neq i}} \sum_{p \in P}\widehat{x}_{ijp}((\sum_{k \in K_{p}}(-c_{p}+\alpha_{jpk}-h_p)s_{ipk})-\beta_{ijp})   \\
					\label{multipliers_cons1}\mbox{s.t.} & & \sum_{k \in K_{p}}(\alpha_{jpk}-c_{p}-h_p)s_{ipk}-\beta_{ijp} \leq \sum_{k \in K_{p}} (\alpha_{ipk}-h_p)s_{ipk}, , \notag \\ 
					& & \qquad \qquad \qquad \qquad \qquad \qquad \qquad \qquad \mbox{ for all } i\neq j \in I, p \in P, k \in K_{p} \mbox{ if } \widehat{x}_{ijp}=1, i= j,  \\
					\label{multipliers_cons2} & & \sum_{k \in K_{p}}(\alpha_{jpk}-c_{p}-h_p)s_{ipk}-\beta_{ijp} \leq \sum_{k \in K_{p}}(\alpha_{j^{*}pk}-c_{o}-h_p)s_{ipk}-\beta_{ij^{*}p} , \notag \\ 
					& & \qquad \qquad \qquad \qquad \qquad \qquad \qquad \qquad \mbox{for all } i,j \in I, p \in P, k \in K_{p}\mbox{ if } \widehat{x}_{ijp}=1, i\neq j, j^{*},\quad \quad \\
					\label{multipliers_cons3} & & \sum_{k \in K_{p}} (\alpha_{ipk}-h_p)s_{ipk} \leq \sum_{k \in K_{p}}(\alpha_{j^{*}pk}-c_{p}-h_p)s_{ipk}-\beta_{ij^{*}p},\notag \\
					& &  \qquad \qquad \qquad \qquad \qquad \qquad \qquad \mbox{ for all } i,j \in I, p \in P, k \in K_{p}\mbox{ if } \widehat{x}_{ijp}=1, i\neq j \neq j^{*},  \\
					\label{multipliers_profit} & & \alpha_{ipk} \leq  r_{p}+h_p, \mbox{ for all } i \in I, p \in P \mbox{ and } k \in K_{p}, \\
					\label{multipliers_cost} & & \beta_{ijp} \leq \sum_{k \in K_{p}}((-c_{p}+r_p-h_p)s_{ipk}), \mbox{ for all } i,j \in I \mbox{ and } p \in P,\\
					\label{multipliers_nonnegative_alpha} & & \alpha_{ipk} \geq  0, \mbox{ for all } i \in I, p \in P \mbox{ and } k \in K_{p}, \\
					\label{multipliers_nonnegative_beta} & & \beta_{ijp} \geq 0, \mbox{ for all } i \in I, j \in I,\mbox{ and } p \in P.
				\end{eqnarray}
			\end{subequations}
			
			\indent The objective function (\ref{multipliers_objective}) is the objective function of the dual problem. Since the solutions $\widehat{\text{x}}$, $\widehat{\text{y}}$, and $\widehat{\text{z}}$ are known, Constraints (\ref{multipliers_cons1})-(\ref{multipliers_nonnegative_beta}) are added to modify the Lagrangian multipliers while retaining the primal solutions. Constraints (\ref{multipliers_cons1}) ensure that if a product $p$ is sent from store $i$ to any other store $j$, then the coefficient of $\widehat{x}_{ijp}$ in the objective function must be less than the coefficient of $\widehat{x}_{iip}$. Similarly, the coefficient of $\widehat{x}_{ijp}$ must be less than the coefficient of any other $\widehat{x}_{ij^{*}p}$, which is guaranteed by Constraints (\ref{multipliers_cons2}). On the other hand, if $\widehat{x}_{iip}=1$, that is, product $p$ remains at its original location, then $\widehat{x}_{ijp}=0$, which is ensured by Constraints (\ref{multipliers_cons3}). Constraints (\ref{multipliers_profit}) guarantee that multipliers $\alpha_{ipk}$ are not greater than corresponding unit revenues plus holding cost to prevent $\widehat{z}_{ipk}$ to be zero. Likewise, Constraints (\ref{multipliers_cost}) set upper bounds on $\beta_{ijp}$. Finally, multipliers $\alpha_{ipk}$ and  $\beta_{ijp}$ must be non-negative which are ensured by Constraints (\ref{multipliers_nonnegative_alpha}) and (\ref{multipliers_nonnegative_beta}), respectively.\\
			
			The optimal solution to this problem is a tighter upper bound as compared to the solution obtained from the relaxed problem. Naturally, this solution provides an upper bound to the optimal solution of the original problem as well.
			
			\subsection{Obtaining lower bounds}

			As mentioned earlier, we obtain lower bounds, i.e., feasible solutions, via a construction heuristic followed by an improvement metaheuristic. The construction heuristic consists of two steps, in the first of which we iteratively connect stores until there is no improvement. We start by dividing all store-product combinations into two groups as sender and receiver based on their stock and demand levels without considering the sizes. For each product, if the stock level in a store is more than its demand, the store is classified as sender; otherwise, it is classified as a receiver. Note that, a store can be either in the sender group or in the receiver group for a product (or, in none of the groups in case the stock and demand levels are equal). We then sequentially connect senders to receivers by selecting products randomly. For each store in the sender group, we find a candidate store from the receiver group that creates the highest profit, i.e., revenue less implied costs. A transfer decision is made if Constraints (\ref{basic_knapsack}) and (\ref{basic_shipment_number2}) remain feasible. After all products are selected, we update the sender and receiver groups considering the current transfers. That is, a store that was initially in the sender group and sends its entire inventory to another store may be included in the receiver group in the next iteration. Moreover, a store in the receiver group can continue to stay in the same group, if it still has needs. Otherwise, it will not be considered as a sender or a receiver. This procedure is repeated until there is no improvement in the solution. In the second step, we further investigate profitable transfers that were not made in the previous step due to Constraints (\ref{basic_shipment_number2}). Now, we search for beneficial transfers by choosing among the destinations that a store is already connected, so that the constraint remains feasible, while the solution is improved. \\
			
			At the improvement stage, we have developed a simulated annealing based metaheuristic. The proposed metaheuristic essentially destroys the current feasible solution by removing a transfer and then repairing it by inserting another transfer. It removes transfers according to three rules that are applied randomly. In the first rule, the transfer to be removed is also selected randomly. The other two rules use the ``residual demand" information for each store-product pair, i.e.,  the difference between the demand and the transfer it receives in the current solution. That is, those with the negative residual demand are the ones that receive more than their demand. The second rule randomly chooses a store-product pair among those that have negative residual demands. And finally, for the third rule we first list all store-product pairs that have negative residual demand. Then we select the product that appears the most in the list and then choose from the stores that also appears the most in the list and paired with this product. After selecting a store-product pair, that transfer is removed and another transfer is inserted while maintaining the feasibility of Constraints (\ref{basic_assignment1}). The destination store is chosen randomly among the ones that have \emph{positive} residual demand. The purpose of these rules is to enable moving to worse as well as better solutions than the current one.\\
			
			The algorithm allows non-improving moves to include diversity as in the simulated annealing (SA) approach. It is adopted as follows: If the profit of the new transfer is greater than or equal to the profit of the removed one, the transfer is accepted. Otherwise, we accept it with probability $ e^{\frac{-(currentProfit - newProfit)}{temperature}}$, where $currentProfit$ and $newProfit$ denote the profits of the removed transfer and the newly added one, respectively and $temperature$ is the current temperature, which is a parameter of SA. Initially, $temperature$ is equal to the total profit of	 current solution so that the $probability$ becomes high and the chance of accepting a worse solution is high. The $temperature$ is decreased at each iteration using the formula $temperature = temperature \times \tau$, where $0 < \tau < 1$ is the $cooling$ $rate$. A $counter$ keeps the number of worse solutions accepted. The $cooling$ $	rate$ is calculated by $cooling$ $ rate = 1/ counter$. The best solution is kept and updated whenever a better solution is found. To avoid being trapped in local optima, the algorithm continues to search from either the best solution or second best solution if there is no improvement in a predefined number of iterations. The algorithm stops if either the total number of iterations reaches to its upper bound or the time limit hast been reached. 
			
			\section{Computational results}\label{CR}
			
			In this section, we report on our computational experiments that consist of two main parts. In the first part, our purpose is to compare the effectiveness of our solution method to that of a commercial solver.  Towards this purpose, the generated instances that are first solved by Gurobi 8.1 and then by our algorithm, which is implemented in Python 3.6 and the subproblems are also solved with Gurobi. All experiments are conducted on a High-Performance Computing (HPC) cluster with Linux system, 44 GB RAM, and 2.40 Ghz processors with two cores. In the second part, our purpose is to develop insights into the effects of the operational constraints on the system performance. Towards this end, we develop another set of instances, all of which are solved by Gurobi.\\
			
			LC Waikiki has about 450 stores (excluding the outlet ones) and about 2,000 items considered for transfer at any week. Unfortunately, problems of this scale cannot be directly handled by Gurobi. Therefore, we have targeted 50 and 100 as the number of stores and 100, 200, 500, and 1000 as the number of products. The number of sizes varies according to the product. However, most of the products have five to 10 different sizes (e.g., S, M, L, XL, and XXL or 28, 30, 32, 34, 36, 38, 40, 42, 44, and 46 for two different products). Thus, in the test problems, the number of sizes is set to either five or 10. Detailed information on the combination of sizes of the instances are given in the first columns of Tables \ref{low}, \ref{medium}, and \ref{high}. In total, we have solved  instances of 14 different size combinations. \\
			
			We randomly generated demands ($d_{ipk}$) and initial stock levels ($s_{ipk}$) from a discrete uniform distribution that is defined between 0 and 10. The prices of the products are set between 20 and 50 Turkish Lira (TL) while the transfer costs for these products are set between 0.4 to 1.5 TL. We generated the selling prices and transfer costs randomly from uniform distributions with the bounds given above. To set the holding cost rate we should have also drawn unit costs but in the interest of simplicity we used the unit revenues. The holding cost per week is taken as the 0.5\% of the unit revenue, which corresponds to about 30\% or less, annually. Although we have not used any real data from LC Waikiki to develop these instances, we have decided on these values upon our conversations with the group that deals with the transfers; hence, we believe our instances are quite realistic.  \\
			
			In order to find the parameters of two of the operational constraints, first we solved each instance  by ignoring all three restrictions. We then found the number of items each store sends and the number of stores it is connected. These are essentially, the maximum values  when there are no operational constraints. Based on these numbers, we then set three levels of $A_i$ and $B_i$ for each store as low (1/3 of the maximum), medium (1/2 of the maximum), and high (2/3 of the maximum).\\
			
\begin{figure}[!b]
				\centering
				\hspace*{-1em}
				\includegraphics[width=0.78\textwidth]{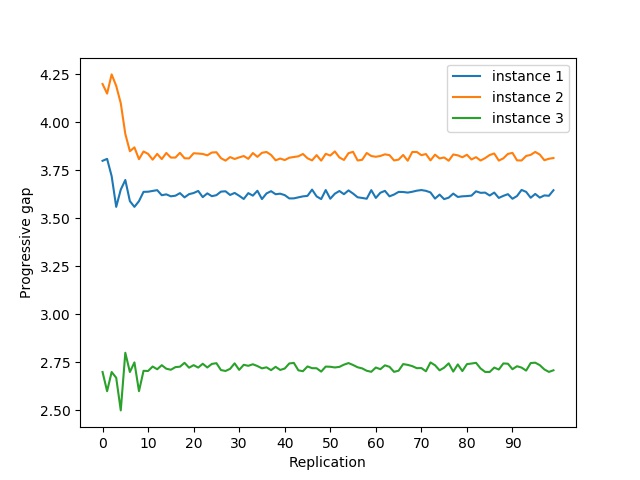}
				\caption{Sensitivity of progressive optimality gap to the number of replications}
				\label{replication}
			\end{figure}

			For each size and ($A_i$, $B_i$) combinations, we have randomly generated 10 instances. We have reached at this number through a small numerical experiment. We took three instance sizes, generated 100 random instances of each, solved them with our algorithm, and computed the average gaps progressively. It turned out that after 10th replication the progressive average of gaps becomes nearly constant for all three sets of instances as illustrated in Figure \ref{replication}. Therefore, we concluded that 10 instances would be enough to have a reliable performance metric in terms of average gaps. As a result, we have generated and tested a total of 14x3x10 = 420 instances.\\
			
			The results are illustrated in Tables \ref{low}, \ref{medium}, and \ref{high}. The first column depicts the size of each instance with respect to the number of stores, products, and sizes. The Gurobi column reports the average of best feasible solutions, the minimum, maximum, and average optimality gaps of 10 replications that Gurobi achieved and the time limit that we set. We have given a one-hour time limit for smaller sized problems, two-hour time limit for medium sized ones, and six-hour time limit to the larger sized instances, in addition to the problem loading times to Gurobi. The metaheuristic column also reports the average of best feasible solutions, the minimum, maximum, and average optimality gaps of 10 replications and the time our algorithms spent to find the upper and lower bounds. We calculated the optimality gaps as $\frac{UB-LB}{LB}$ where $UB$ and $LB$ represent the upper and lower bounds found by Gurobi and our method.\\

			\begin{table}[t]
				\centering
				\scriptsize
				\caption{Results for Low level of $A_i$ and $B_i$}
				\hspace*{-3em}
				\begin{tabular}{ccccccccccc}
					\hline
					& \multicolumn{5}{c}{Gurobi}            & \multicolumn{5}{c}{Metaheuristic} \bigstrut[b]\\
					\cline{2-11}    Instance & Average & \multicolumn{3}{c}{Gap} & Average  & Average  & \multicolumn{3}{c}{Gap} & Average  \bigstrut[t]\\
					$I$-$P$-$K$ &  Lower Bound & Min   & Average & Max   & Runtime & Lower Bound & Min   & Average & Max   & Runtime \bigstrut[b]\\
					\hline
					50-100-5 & 848,013 & 0.71  & 0.79  & 0.89  & 3600 & 841,322 & 1.04  & 2.07  & 3.67  &  1285 \bigstrut[t]\\
					50-100-10 & 307,259 & 0.00  & 0.00  & 0.01  & 7.04  & 302,956 & 0.50  & 2.41  & 5.33  & 1017\\
					50-200-5 & 1,519,767 & 0.78  & 0.87  & 1.01  & 3600& 1,500,891 & 1.06  & 3.00  & 5.38  & 1424 \\
					50-200-10 & 935,098 & 0.79  & 1.09  & 1.37  & 3600 & 924,970 & 1.77  & 3.39  & 8.32  & 1318 \\
					50-500-5 & 2,317,320 & 282.10 & 288.27 & 296.74 & 3600 & 6,244,111 & 4.22  & 4.48  & 4.76  & 2001 \\
					50-500-10 & 2,417,042 & 0.34  & 0.40  & 0.45  & 3600 & 2,398,677 & 0.98  & 2.40  & 3.65  & 1960 \\
					50-1000-5 & 4,663,613 & 262.45 & 284.38 & 302.65 & 3600 & 11,074,571 & 2.50  & 3.76  & 4.29  & 1966 \\
					50-1000-10 & 3,507,723 & 0.09  & 0.14  & 0.22  & 3600 & 3,487,356 & 0.76  & 1.04  & 1.44  & 2082 \\
					100-100-5 & 950,596 & 283.04 & 287.71 & 293.21 & 7200& 2,695,145 & 3.47  & 5.30  & 6.20  & 1492 \\
					100-100-10 & 1,890,957 & 0.22  & 0.31  & 0.41  & 7200 & 1,882,898 & 0.57  & 1.58  & 3.29  & 3707 \\
					100-500-5 & 4,666,991 & 284.15 & 287.76 & 294.74 & 7200 & 12,141,278 & 6.10  & 6.57  & 6.96  & 4313 \\
					100-500-10 & 7,809,377 & 0.71  & 0.91  & 1.11  & 7200 & 7,714,316 & 0.85  & 3.52  & 5.16  & 4414 \\
					100-1000-5 & 9,303,094 & 284.71 & 288.34 & 291.90 & 21600 & 18,948,867 & 0.49  & 0.94  & 1.24  & 9245 \\
					100-1000-10 & 12,997,254 & 0.39  & 0.99  & 1.28  & 21600 & 12,906,313 & 1.33  & 2.45  & 3.22  & 6475 \\
					\hline
				\end{tabular}%
				\label{low}%
			\end{table}%

			\begin{table}[!b]
				
				\centering
				\scriptsize
				\caption{Results for Medium level of $A_i$ and $B_i$}
				\hspace*{-3em}
				\begin{tabular}{ccccccccccc}
					\hline
					& \multicolumn{5}{c}{Gurobi}            & \multicolumn{5}{c}{Metaheuristic} \bigstrut[b]\\
					\cline{2-11}    Instance & Average & \multicolumn{3}{c}{Gap} & Average  & Average  & \multicolumn{3}{c}{Gap} & Average  \bigstrut[t]\\
					$I$-$P$-$K$ &  Lower Bound & Min   & Average & Max   & Runtime & Lower Bound & Min   & Average & Max   & Runtime \bigstrut[b]\\
					\hline
					50-100-5 & 1,023,543 & 0.42  & 0.49  & 0.56  & 3600 & 1,015,396 & 0.82  & 2.57  & 4.95  & 1274 \bigstrut[t]\\
					50-100-10 & 547,922 & 0.01  & 0.01  & 0.01  & 169 & 542,032 & 0.63  & 1.74  & 3.05  & 1181 \\
					50-200-5 & 1,794,245 & 0.72  & 0.83  & 0.95  & 3600 & 1,760,781 & 1.60  & 3.58  & 5.37  & 1404 \\
					50-200-10 & 1,249,620 & 0.47  & 0.56  & 0.66  & 3600 & 1,224,849 & 1.09  & 3.36  & 6.11  & 1366 \\
					50-500-5 & 5,085,082 & 6.65  & 90.41 & 287.53 & 3600 & 6,501,742 & 3.70  & 3.96  & 4.23  & 2000 \\
					50-500-10 & 3,192,027 & 0.08  & 0.10  & 0.13  & 3600 & 3,149,863 & 0.73  & 3.53  & 4.67  & 1982 \\
					50-1000-5 & 4,663,613 & 283.04 & 287.71 & 293.21  & 3600 & 12,553,494 & 1.12  & 1.80  & 2.93  & 1974 \\  
					50-1000-10 & 6,070,388 & 0.10  & 0.14  & 0.19  & 3600 & 5,984,386 & 1.08  & 1.75  & 2.56  & 2075 \\
					100-100-5 & 3,008,779 & 0.52  & 0.82  & 1.18 & 7200 & 2,997,312 & 2.57  & 4.49  & 5.21  & 1475 \\
					100-100-10 & 1,890,938 & 0.07  & 0.08  & 0.10  & 7200 & 1,879,030 & 0.34  & 2.79  & 4.13  & 3774 \\
					100-500-5 & 4,636,738 & 284.96 & 289.79 & 294.74 & 7200 & 12,713,377 & 4.81  & 5.30  & 5.76  & 4754 \\
					100-500-10 & 9,447,082 & 0.94  & 1.06  & 1.17  & 7200 & 9,363,384 & 2.09  & 2.74  & 3.43  & 4514 \\
					100-1000-5 & 9,325,113 & 284.71 & 287.41 & 289.97 & 21600 & 20,854,343 & 1.83  & 2.50  & 2.92  & 9028 \\
					100-1000-10 & 15,871,661 & 0.01  & 0.42  & 1.22  & 21600 & 15,780,238 & 1.11  & 1.55  & 2.52  & 6109 \\ \hline
				\end{tabular}%
				\label{medium}%
			\end{table}%

			\begin{table}[t]
				\centering
				\scriptsize
				\caption{Results for High level of $A_i$ and $B_i$}
				\hspace*{-3em}
				\begin{tabular}{ccccccccccc}
					\hline
					& \multicolumn{5}{c}{Gurobi}            & \multicolumn{5}{c}{Metaheuristic} \bigstrut[b]\\
					\cline{2-11}    Instance & Average & \multicolumn{3}{c}{Gap} & Average  & Average  & \multicolumn{3}{c}{Gap} & Average  \bigstrut[t]\\
					$I$-$P$-$K$ &  Lower Bound & Min   & Average & Max   & Runtime & Lower Bound & Min   & Average & Max   & Runtime \bigstrut[b]\\
					\hline
					50-100-5 & 1,171,178 & 0.33  & 0.39  & 0.48  & 15 & 1,154,752 & 1.79  & 2.96  & 5.17  & 1290 \bigstrut[t]\\
					50-100-10 & 819,047 & 0.01  & 0.02  & 0.06  & 2845 & 805,510 & 0.37  & 2.13  & 3.10  & 1333 \\
					50-200-5 & 2,042,125 & 0.58  & 0.79  & 1.37  & 3600 & 2,017,624 & 1.52  & 3.36  & 5.10  & 1448 \\
					50-200-10 & 1,507,772 & 0.31  & 0.37  & 0.46  & 3600 & 1,496,739 & 1.31  & 2.63  & 3.88  & 1423 \\
					50-500-5 & 7,151,269 & 0.12  & 0.28  & 0.41  & 3600 & 7,120,120 & 0.72  & 1.06  & 1.97  & 1778 \\
					50-500-10 & 4,456,310 & 0.01  & 0.01  & 0.01  & 154 & 4,296,653 & 3.37  & 4.28  & 4.77  & 2131 \\
					50-1000-5 & 12,891,943 & 0.07  & 0.12  & 0.16  & 3600 & 12,839,782 & 0.46  & 0.75  & 1.13  & 1971 \\
					50-1000-10 & 8,840,304 & 0.01  & 0.01  & 0.02  & 1963 & 8,462,620 & 1.46  & 2.12  & 2.98  & 2064 \\
					100-100-5 & 3,021,758 & 0.20  & 0.30  & 0.46  & 7200& 3,010,490 & 3.76  & 4.55  & 5.64  & 1464 \\
					100-100-10 & 1,891,297 & 0.53  & 0.56  & 0.63  & 7200 & 1,872,998 & 4.41  & 5.23  & 5.96  & 3647 \\
					100-500-5 & 14,830,650 & 0.49  & 1.27  & 2.16  & 7200 & 14,685,243 & 2.66  & 3.87  & 5.58  & 4201 \\
					100-500-10 & 9,492,130 & 0.02  & 0.14  & 0.29  & 7200 & 9,373,087 & 1.10  & 2.19  & 3.02  & 4512 \\
					100-1000-5 & 29,909,041 & 0.56  & 0.82  & 1.56  & 21600 & 29,313,762 & 1.55  & 3.39  & 5.72  & 9014 \\
					100-1000-10 & 19,301,509 & 0.00  & 0.01  & 0.01  & 11479 & 19,217,184 & 0.65  & 1.09  & 1.43  & 6004 \\ \hline
				\end{tabular}%
				\label{high}%
			\end{table}%
			
			The results show that our algorithm is comparable to, and in some cases much more effective than, Gurobi. First of all, our algorithm spends about one-third to one-half of the time that we give to Gurobi excluding the time it takes to load the problem to Gurobi, which could be rather substantial in larger instances. In terms of the optimality gaps, the results are somewhat mixed. There are many instance sets, for which Gurobi found better solutions (lower bounds) than our method did. However, Gurobi's solutions deteriorate faster than our method with increasing problem size. Although our approach also suffers, it can solve most medium-sized problems with around 1\% optimality gaps and large-sized problems with a maximum gap of about 7\%.\\
			
			The most important problem with Gurobi, however, is that it is rather unreliable. In some cases the solutions it found were just terrible, with optimality gaps hovering around 300\%.  Upon close inspection, we noticed that these poor results belong to the instances where there are five different sizes, whereas the instances with the same number of stores and products but 10 different sizes for each products, the behavior was quite the opposite. This was rather puzzling; after all, the latter instances are of larger size, but with closer inspection we were able to conclude that it was the combination of several factors that led to this unexpected results. First of all, since the transportation costs are quite low as compared to the revenues, as long as it is revenue-improving the optimal solution tends to have large number of transfers. Secondly, when there are 10 sizes for each product, there are fewer profitable opportunities for transfers as compared to the same number of stores and products with five sizes for each product. This might seem unclear at first, but single-destination constraint is mainly responsible for these results. For example, if there is only one size, then there would be many profitable opportunities for transfer. When there are two sizes, the opportunities would diminish because there would be more sales opportunities at the original sources and there would be fewer alternative stores that would have demand for both sizes. This would be even more prominent with increasing the number of sizes. This is indeed what we have observed in a simple experiment that we have conducted with 20 stores, 100 products and no operational constraints. We then set the number of sizes as 2, 5, 8, 10, and 15 and have randomly drawn 10 instances for each size. Table \ref{siz} reports the results, which confirm our intuition. Hence, since there are far fewer transfers in the optimal solution as the number of sizes increases most stores do not perform any transfers but keep at the source. As a result, Gurobi can eliminate a substantial number of potential transfers across the stores and finds solutions easier in instances where there are 10 sizes as compared to five. To conclude, while our approach does not produce better solutions than Gurobi all the time, due to its robustness to problem characteristics, it is a much better alternative of the two.\\
			
			\begin{table}[!t]
				\centering
				\caption{Effect of the number of sizes on the number of transfers}
				\hspace*{-1em}
				\scalebox{0.75}{
					\begin{tabular}{ccccccccccc}
						\hline
						$K$     & Instance 1 & Instance 2 & Instance 3 & Instance 4 & Instance 5 & Instance 6 & Instance 7 & Instance 8 & Instance 9 & Instance 10 \bigstrut\\
						\hline
						2     & 1,038  & 1,074  & 1,009  & 1,013  & 1,031  & 1,030  & 1,069  & 1,011  & 998   & 1,038 \bigstrut[t]\\
						5     & 616   & 626   & 624   & 650   & 630   & 625   & 634   & 653   & 612   & 621 \\
						8     & 327   & 371   & 334   & 366   & 352   & 353   & 315   & 335   & 315   & 354 \\
						10    & 208   & 226   & 220   & 220   & 206   & 231   & 204   & 244   & 207   & 237 \\
						15    & 42    & 55    & 65    & 57    & 48    & 53    & 59    & 57    & 78    & 50 \bigstrut[b]\\
						\hline
					\end{tabular}}
					\label{siz}%
				\end{table}%

\begin{figure}[!h]
					\hspace*{-6em}
					\begin{minipage}{0.5\linewidth}
						\includegraphics[width = 1.3\textwidth]{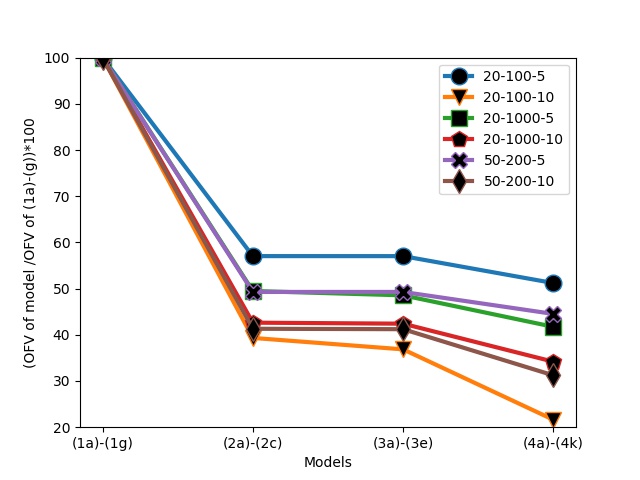}
					\end{minipage}
					\hfill
					\hspace*{-6em}
					\begin{minipage}{0.505\linewidth}
						\includegraphics[width = 1.3\textwidth]{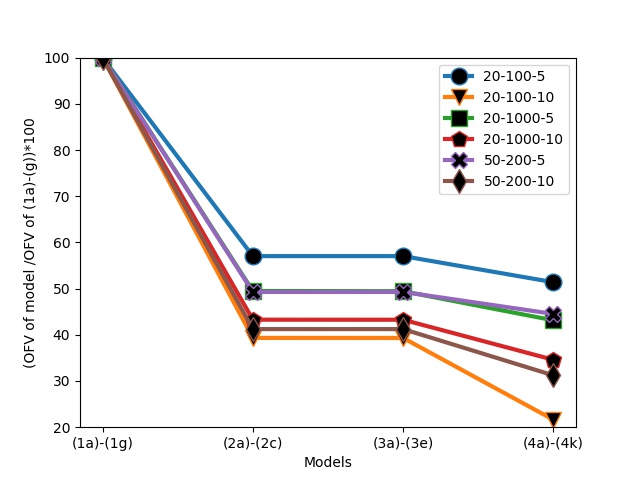}
					\end{minipage}
					\hspace*{-6em}
					\begin{minipage}{0.5\linewidth}
						\includegraphics[width = 1.3\textwidth]{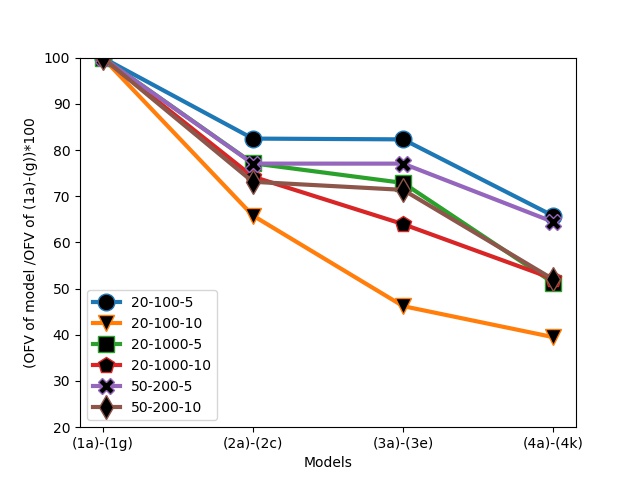}
					\end{minipage}
					\hfill
					\hspace*{-6em}
					\begin{minipage}{0.505\linewidth}
						\includegraphics[width = 1.3\textwidth]{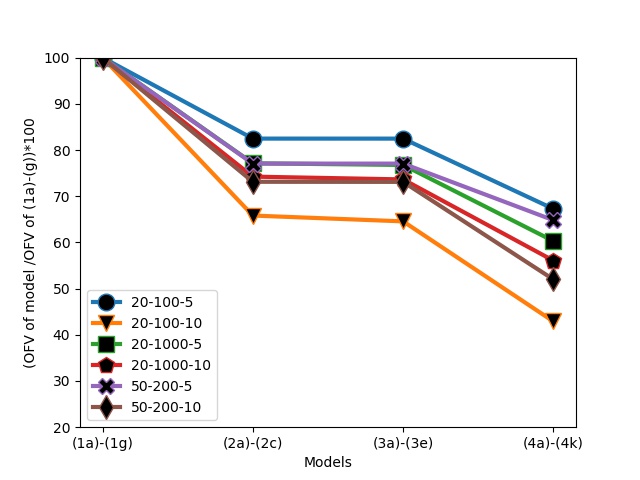}
					\end{minipage}
					\caption{Effect of adding capacity and single destination constraints on the objective function value (OFV). Top-left: low $A_i$ and low $B_i$, top-right: low $A_i$ and high $B_i$, bottom-left: high $A_i$ and low $B_i$, and bottom-right: high $A_i$ and high $B_i$.}
					\label{managerial}
				\end{figure}

				In the second part of our experiments, our purpose is to shed some light into the effects of particular operational restrictions used by LC Waikiki. The restrictions include the two capacity constraints and the single-destination policy. Towards this end, we considered six problem sizes as illustrated in Figure \ref{managerial}. We randomly generated 10 instances for each of these problem sizes. As we have done in the first part, we have created further instances based on how tight the capacity constraints are. Similarly, we first solved each instance  by ignoring all three restrictions, found the maximum values $A_i$ and $B_i$ can take for each store, and then created four combinations of ($A_i$, $B_i$) by setting them to either ``low" (1/3 of the maximum) or ``high" (2/3 of the maximum). Therefore, altogether we have solved a total of 240 instances. \\
				
				The unconstrained version of the problem, i.e., \eqref{relaxed_objective}-\eqref{relaxed_binary_z}, is solved first, and then the problem with the transfer capacity constraint, i.e., \eqref{relaxed_objective_ex1}-\eqref{all_const_ex1}, which is followed by the problem with both capacity constraints, i.e., \eqref{relaxed_objective_ex2}-\eqref{all_const_ex2}, and finally, the full problem \eqref{basic_objective_nl}-\eqref{basic_nonzero_z} is solved. All these problems are solved with Gurobi. Although not all problems are solved to optimality, the gaps are rather small, so the results are quite reliable. \\
				
				Figure \ref{managerial} depicts the summary results. Each graph in the figure contains results based on the combinations of ($A_i$, $B_i$) pairs. In the figure, optimum solution of each unconstrained instance is normalized to 100 and the objective functions of each instance's constrained versions are found as percentage of the optimal value of the unconstrained version. The graphs report the averages of these percentages over 10 instances. The effect of the constraint on the number of SKUs that can be transferred is quite clear. The optimal values reduces to roughly 40-60\% and 65-85\% of the maximum possible for low and high $A_i$'s, respectively. The addition of the second constraint that restricts the number of stores has almost no further deteriorating effect. Hence, it appears that the first constraint is already restrictive enough. The only exception to these results is when $A_i$ is high \emph{and} there are 10 sizes for each product. As the bottom-left graph shows, although adding the second constraint has a very slight effect when there are five sizes, its effect is substantial when it is 10. This result is actually quite intuitive because when there are 10 sizes there are simply more opportunities to match demand and supply as there are more sizes (because there is no single-destination constraint yet). Therefore, restricting the number of stores greatly eliminates those opportunities. \\
				
				We can observe that the single destination constraint has a substantial deteriorating effect on the profit. Depending on the cases, it has roughly an additional 5-25\% negative effect on the maximum profit. The effect is somewhat less when $A_i$ is low, presumably the first constraint already has a great effect on reducing the number of destinations. However, the effect is quite substantial when $A_i$ is high regardless of the value of $B_i$.\\
				
				To summarize the managerial implications, we can conclude that in general, i) it is the restriction on the total number of SKUs that can be transferred rather  than the restrictions on the number of destinations that has a more negative effect, ii) single destination constraint has a more detrimental effect when the capacities are less restrictive, and finally, iii) all the negative effects are usually more pronounced when there are larger number of sizes per product. Therefore, if a firm wishes to relax the single-destination constraint, it should start from products with larger number of sizes and accompany it relaxing restrictions on the transfer capacities.

				\section{Concluding remarks}\label{Concl}
				In this paper we introduce a novel proactive transshipment problem motivated by the practice at the largest fast fashion retailer in Turkey, LC Waikiki. The company, after allocating the initial inventory to over 450 stores and observing sales for a few weeks, engages in lateral transshipments among the stores. When a product has different sales performances across the stores, lateral transshipments can improve the overall system performance. Not only such a practice helps the company better match supply with demand but also eliminates additional handling and transportation operations at its central depot.  \\

				Its large scale and particular operational restrictions necessitate the development of a novel model. We formulate the transfer problem of LC Waikiki as a mixed integer linear programming problem. With around 450 stores, 2,000 products, and a variety of operational constraints, this problem becomes a very large mixed integer program and solving it optimally becomes a challenge. Therefore, we have developed a simulated annealing based metaheuristic to solve the problem. We also applied Lagrangian relaxation with a primal-dual approach to obtain sharp upper bounds on the optimal solution of the original problem. We generated 420 problem instances of varying sizes to evaluate the performance of the proposed algorithm against the commercial optimizer Gurobi. Each instance is solved by the proposed algorithm and Gurobi. The results show that although the solutions prescribed by Gurobi are better than ours in small-size instances and those with loose capacities, the proposed algorithm outperforms Gurobi in instances that are characterized by having a large number of potentially beneficial transfers and tighter capacity constraints. These instances are the ones where the combinatorial nature of the problem becomes the most challenging. Gurobi fails spectacularly in these instances, while our algorithm performs without a significant loss in its performance. Hence, our algorithm is quite robust to changing program characteristics. Finally, our algorithm is also quicker in finding solutions, spending only about one-third to one-half of the time spent by Gurobi. This feature makes our approach particularly attractive when companies need speedy solutions to these problems.  \\
				
				We have also conducted a carefully designed numerical experiment to uncover the effect of the particular operational constraints of the company. First, we have solved the instances without any of those constraints and found the maximum potential gross revenue (i.e., the base gross revenue) that can be obtained with transshipments. We have then added those constraints one by one to observe their effects. We observe that constraints on the total number of products a store can send has a significant impact and depending on how tight those restrictions are, may reduce the gross revenue to around 40-85\% of its base level. The second capacity restriction, i.e., the number stores that a store can make shipments, usually has very little negative impact after the first capacity restriction is already imposed. Finally, we have also measured the effect of the single-destination practice and have found that it may reduce the revenues by another 5-25\% of the base revenue depending on how tight the capacity constraints are. When the capacity constraints are already tight, the negative impact of this practice is quite small, but where the capacities are loose the negative impact of this practice is quite significant. We have also observed that the number of sizes also plays a significant role in these results. In our sample instances we have used five and 10 as the number of sizes. Naturally, when there are 10 sizes of products, single-destination practice renders much fewer number of transfers as potentially beneficial and therefore, this practice becomes more detrimental to the base revenue when the number of sizes increases.\\

				We had to make a number of simplifying assumptions to effectively deal with this very large problem that has complicating operational constraints. Therefore, there are a number of avenues for future research. While one may have accurate demand forecasts as in this case, there are always forecast errors, and therefore, considering demand uncertainty is naturally an important extension to this study. In a similar vein, initial shipment decisions under demand uncertainty may also be  considered jointly with the transfer decisions. Another potentially important avenue is to develop integrated models that include transfer decisions as well as markdown decisions, an avenue that we are currently pursuing. Finally, the frequency at which the collections are renewed can also be made jointly with transfer as well as other logistical decisions. \\
				
				\section*{Acknowledgments}
				The authors first thank to their industrial partner, LC Waikiki and their key collaborator, Allocation Strategy and Process Development team- in particular, Murat Aksu (the team's previous director), Serkan Ulukaya, Mustafa Tilkat, Yavuz Selim Dola\c{s}, and Muhammed Can Konur. They are also particularly grateful to other employees at LC Waikiki, including Elif Hilal Akar, Aytu\u{g} Alten, and Sinan Erg\"{u}n for their support.
				
				\section*{References}

			\end{document}